\documentclass{article}

\usepackage{graphicx}

\def\be{\begin{equation}}
\def\ee{\end{equation}}
\def\ff#1{\mbox{\boldmath $#1$} }
\def\a{\alpha}
\def\b{\beta}
\def\lam{\lambda}
\def\e{\epsilon}
\def\x{\mbox{\boldmath $x$} }
\newcommand{\mat}[1]{{\left( \begin{array}{cccc}#1\end{array}\!\right)}}
\newcommand{\brak}[1]{\left\{  \begin{array}{lllll} #1 \end{array} \right. }

\def\={\approx}
\newcommand{\kkk}[1]{$^{{\scriptsize~\cite{#1}}}$}

\begin{document}

\title{\bf Efficiency Analysis of Swarm Intelligence and Randomization Techniques}

\author{Xin-She Yang \\
Mathematics and Scientific Computing, \\
National Physical Laboratory,   \\
Teddington, Middlesex TW11 0LW, UK.  }

\date{}

\maketitle

\begin{abstract}

Swarm intelligence has becoming a powerful technique in  solving design
and scheduling  tasks. Metaheuristic algorithms are an integrated part of this paradigm,
and particle swarm optimization is often viewed as an important landmark.
The outstanding performance and efficiency of swarm-based algorithms
inspired many new developments, though mathematical understanding of
metaheuristics remains partly a mystery. In contrast to the classic
deterministic algorithms, metaheuristics such as PSO always use some
form of randomness, and such randomization now employs various techniques.
This paper intends to review and analyze some of the convergence and efficiency
associated with metaheuristics such as firefly algorithm, random walks,
and L\'evy flights. We will discuss how these techniques are used and
their implications for further research. \\

{\bf Citation details:}
X. S. Yang, Efficiency analysis of swarm intelligence and randomization techniques,
{\it Journal of Computational and Theoretical Nanoscience}, Vol. 9, No. 2, pp. 189--198 (2012).

\end{abstract}

\section{Swarm Intelligence}

The emerging trend of cognitive informatics is to use swarm intelligence
to find solutions to challenging design problems and real-world applications.
A representative example is the development of particle swarm optimization(PSO),
which can be considered as a classic landmark in using swarm intelligence
or collective intelligence\kkk{Kennedy}. Another excellent example is the
ant colony optimization which mimics the collective behaviour of social insects\kkk{Dorigo}.
These algorithms are metaheuristic and swarm-based, though they
may have some similarities with genetic algorithms\kkk{Bell,Kennedy,Enge},
but it is much simpler because they do not use mutation/crossover
operators. Genetic algorithms are population-based, which was based on
the evolution theory and survival of the fittest\kkk{Holland}.
Similarly, PSO is also population-based, but it was inspired by the
swarm behaviour of fish and birds, and therefore falls into the different
category from genetic algorithms.

Swarm-based or population-based algorithms form the majority of nature-inspired
metaheuristic algorithms\kkk{Blum,Talbi,Yang,Yang4,Yang5} with
a wider range of applications\kkk{Ayesh,Gross,Kennedy2}.
Population-based algorithms are not necessarily `intelligent'. In contrast,
the the trajectory-based algorithms such as simulated annealing\kkk{Kirk} may
behave slightly `intelligent' somehow.  Even an
algorithm is called swarm-based, but it is not really `intelligent' in the normal
sense. What we really mean is that the algorithm intends to mimic the
swarm behaviour and work in an idealized manner, in the hope that it can provide
better and more efficient search for finding tough optimization problems,
and thus showing some sort of `intelligence' comparing with the classical
mechanistical algorithms.

Since the appearance of swarm intelligence algorithms such as PSO in the 1990s,
more than a dozen of new metaheuristic algorithms have been developed\kkk{Dorigo,Chow,Geem,Yang,Yang5},
and these algorithms have been applied to almost all areas of optimization, design,
scheduling and planning, data mining, machine intelligence, and many others.
Thousands of research papers and dozens of books have been published\kkk{Talbi,Yang2},

Despite the rapid development of swarm-based metaheuristics, their mathematical
analysis remains partly unsolved, and many open problems need urgent attention.
This difficulty is largely due to the fact the interaction of various components
in metaheuristic algorithms are highly nonlinear, complex, and stochastic.
Studies have attempted to carry out convergence analysis using the theory and methods of
dynamical systems, and some important results concerning PSO were obtained\kkk{Clerc}.
However, for other metaheuristics such as firefly algorithms and ant colony optimization,
it remains an active, challenging topic.
On the other hand, even we have not proved or cannot prove their convergence, we still can
compare their performance of various algorithms. This has indeed formed a majority of
current research in algorithm development in the research community of optimization
and machine intelligence\kkk{Enge,Talbi,Yang2}.

Another important issue is the various randomization techniques used for modern
metaheuristics, from simple randomization such as uniform distribution
to random walks, and to more elaborate L\'evy flights\kkk{Blum,Pav,Yang2}. There
is no unified approach to analyze these mathematically. In this paper, we intend
to review the convergence of swarm-based algorithms including PSO
and firefly algorithms. We also try to analyze the mathematical and
statistical foundation for randomization techniques from simple random walks
to L\'evy flights. As the result of the analysis, we will study the proper
step size for randomization in metaheuristics. Finally, we will discuss the
implication and further research topics in swarm intelligence and their applications.

\section{Particle Swarm Optimization}

\subsection{Standard PSO}

Particle swarm optimization (PSO) was developed by Kennedy and
Eberhart in 1995\kkk{Kennedy,Kennedy2}, based on the swarm behaviour such
as fish and bird schooling in nature. Since then, PSO has
generated much wider interests, and forms an exciting, ever-expanding
research subject, called swarm intelligence. PSO has been applied
to almost every area in optimization, computational intelligence,
and design/scheduling applications.
There are at least two dozens of PSO variants and most of them
were developed by introducing or varying certain components \kkk{Cui,Cui2}.
Hybrid algorithms by combining PSO with other existing algorithms
are also increasingly popular.

In essence, PSO searches the design space of an optimization problem
by adjusting the trajectories of individual agents,
called particles, as the piecewise paths formed by positional
vectors in a quasi-stochastic manner. The movement of a swarming particle
consists of two major components: a stochastic component and a deterministic component.
Each  particle is attracted toward the position of the current global best
$\ff{g}^*$ and its own best location $\x_i^*$ in history,
while at the same time it has a tendency to move randomly.

Let $\x_i$ and $\ff{v}_i$ be the position vector and velocity for
particle $i$, respectively. The new velocity vector is determined by the
following formula
\be \ff{v}_i^{t+1}= \ff{v}_i^t  + \a \ff{\e}_1  \odot
[\ff{g}^*-\x_i^t] + \b \ff{\e}_2 \odot [\x_i^*-\x_i^t].
\label{pso-speed-100}
\ee
where $\ff{\e}_1$ and $\ff{\e}_2$ are two random vectors, and each
entry taking the values between 0 and 1. The Hadamard product of two
matrices $\ff{u \odot v}$ is defined as the entrywise product, that
is $[\ff{u \odot v}]_{ij}=u_{ij} v_{ij}$.
The parameters $\a$ and $\b$ are the learning parameters or
acceleration constants, which can typically be taken as, say, $\a\=\b \=2$.

The initial locations of all particles should distribute relatively
uniformly so that they can sample over most regions, which
is especially important for multimodal problems. The initial velocity
of a particle can be taken as zero, that is, $\ff{v}_i^{t=0}=0$.
The new position can then be updated by
\be \x_i^{t+1}=\x_i^t+\ff{v}_i^{t+1}. \label{pso-speed-140} \ee
Although $\ff{v}_i$ can be any
values, it is usually bounded in some range $[\ff{v}_{\min}, \ff{v}_{\max}]$.

There are many variants which extend the standard PSO
algorithm, and the most noticeable improvement
is probably to use inertia function $\theta
(t)$ so that $\ff{v}_i^t$ is replaced by $\theta(t) \ff{v}_i^t$
\be \ff{v}_i^{t+1}=\theta \ff{v}_i^t  + \a \ff{\e}_1  \odot
[\ff{g}^*-\x_i^t] + \b \ff{\e}_2 \odot [\x_i^*-\x_i^t],
\label{pso-speed-150}
\ee
where $\theta$ takes the values between 0 and 1\kkk{Chat}. In the simplest case,
the inertia function can be taken as a constant, typically $\theta \= 0.5 \sim 0.9$.
This is equivalent to introducing a virtual mass to stabilize the motion
of the particles, and thus the algorithm is expected to converge more quickly.

\subsection{Accelerated PSO}

The standard particle swarm optimization uses both the current global best
$\ff{g}^*$ and the individual best $\ff{x}^*_i$. The reason of using the individual
best is primarily to increase the diversity in the quality solutions, however,
this diversity can be simulated using some randomness. Subsequently, there is
no compelling reason for using the individual best, unless the optimization
problem of interest is highly nonlinear and multimodal.

A simplified version which could accelerate the convergence of the algorithm is to use the global best only\kkk{Yang}. Thus, in the accelerated particle swarm optimization, the
velocity vector is generated by a simpler formula
\be \ff{v}_i^{t+1}=\ff{v}_i^t + \a (\ff{\e}-1/2) + \b (\ff{g}^*-\x_i^t), \label{pso-sys-10} \ee
where $\ff{\e}$ is a random variable with values from 0 to 1. Here the shift
1/2 is purely out of convenience. We can also use a standard normal
distribution  $\a \ff{\e}_n$ where $\ff{\e}_n$ is drawn from $N(0,1)$
to replace the second term.
The update of the position
is simply \be \x_i^{t+1}=\x_i^t + \ff{v}_i^{t+1}.  \label{pso-sys-20} \ee In order to
increase the convergence even further, we can also write the
update of the location  in a single step
\be  \x_i^{t+1}=(1-\b) \x_i^t+\b \ff{g}^* +\a \ff{\e}_n.  \ee
This simpler version will give the same order of convergence.
The typical values for this accelerated PSO are
$\a \=0.1 \sim 0.4$ and $\b \=0.1 \sim 0.7$, though
$\a \=0.2$ and $\b \=0.5$ can be taken as the initial values for most unimodal objective functions.
It is worth pointing out that the parameters $\a$ and $\b$ should in general be related to the scales of the independent variables $\x_i$ and the search domain.

A further improvement to the accelerated PSO is to reduce the randomness
as iterations proceed.
This means   that we can use a monotonically decreasing function such as
\be \a =\a_0 e^{-\gamma t}, \ee
or \be \a=\a_0 \gamma^t, \quad (0<\gamma<1), \ee
where $\a_0 \=0.5 \sim 1$ is the initial value of the randomness parameter.
Here $t$ is the number of iterations or time steps.
$0<\gamma<1$ is a control parameter. This is similar to the geometric cooling
schedule used in simulated annealing\kkk{Kirk}.

As we will see that the convergence of the standard PSO can be
guaranteed under appropriate conditions. The accelerated PSO can
also have guaranteed convergence, which forms the basis
for the convergence of the firefly algorithm algorithm to be introduced later. Interestingly, some chaotic behaviour may emerge under the right conditions, which may be used as advantage
for efficient randomization.

\subsection{Convergence Analysis}

The first convergence analysis of PSO was carried out by Clerc and Kennedy
in 2002\kkk{Clerc} using the theory of dynamical systems.
Mathematically, if we ignore the random factors, we can view the system
formed by (\ref{pso-speed-100}) and (\ref{pso-speed-140}) as a dynamical system.
If we focus on a single particle $i$ and imagine that there is only one particle in this system,
then the global best $\ff{g}^*$ is the same as its current best $\x_i^*$. In this case, we have
\be \ff{v}_i^{t+1} = \ff{v}_i^t + \gamma (\ff{g}^*-\x_i^t), \quad \gamma=\a+\b, \ee
and
\be \x_i^{t+1} =\x_i^t + \ff{v}_i^{t+1}. \ee
Considering the 1D dynamical system for particle swarm optimization,
we can replace $\ff{g}^*$ by a parameter constant $p$ so
that we can see if or not the particle of interest will converge towards $p$.
By setting $u_t=p-x(t+1)$
and using the notations for dynamical systems, we have a simple dynamical system
\be v_{t+1}=v_t + \gamma u_t, \quad u_{t+1} =-v_t + (1-\gamma) u_t, \ee
or
\be Y_{t+1} = A Y_t, \quad A=\mat{ 1 & \gamma \\ -1 & 1-\gamma}, \quad Y_t=\mat{v_t \\ u_t}. \ee
The general solution of this dynamical system can be written as
\be Y_t=Y_0 \exp[A t]. \ee
The main behaviour of this system can be characterized by the eigenvalues $\lam$ of $A$
\be \lam_{1,2}=1-\frac{\gamma}{2} \pm \frac{\sqrt{\gamma^2 - 4 \gamma}}{2}. \ee
It can be seen clearly that $\gamma=4$ leads to a  bifurcation.

Following a straightforward analysis of this dynamical system,
we can have three cases. For $0 < \gamma <4$, cyclic and/or quasi-cyclic
trajectories exist. In this case, when randomness is gradually reduced,
some convergence can be observed.
For $\gamma>4$, non-cyclic behaviour can be expected and
the distance from $Y_t$ to the center $(0,0)$ is monotonically increasing with $t$.
In a special case $\gamma=4$, some convergence behaviour can be observed.
For detailed analysis, please refer to Clerc and Kennedy\kkk{Clerc}.
Since $p$ is linked with the global best, as the iterations continue,
it can be expected that all particles will aggregate towards the
the global best.

\section{Firefly Algorithm}

Firefly Algorithm (FA) was developed by Xin-She Yang at Cambridge University in late 2007
and early 2008\kkk{Yang,Yang3}, which was based on the flashing patterns and behaviour
of fireflies. In essence, FA uses the following three idealized rules:

\begin{itemize}
\item Fireflies are unisex so that one firefly will be attracted to other fireflies
regardless of their sex;

\item The attractiveness is proportional to the brightness and they both
decrease as their distance increases. Thus for any two flashing fireflies, the less brighter one will move towards the brighter one.  If there is
no brighter one than a particular firefly, it will move randomly;

\item The brightness of a firefly is determined by the landscape of the
objective function.
\end{itemize}

For a  maximization problem, the brightness can simply be proportional
to the value of the objective function. Other forms of brightness can be defined in a similar
way to the fitness function in genetic algorithms.

\subsection{Variations of Attractiveness}

As both light intensity and attractiveness affect the movement of fireflies
in the firefly algorithm, we have to define their variations. For simplicity,
we can always assume that the attractiveness of a firefly is
determined by its brightness which in turn is associated with the
encoded objective function. In the simplest case for maximum optimization problems,
the brightness $I$ of a firefly at a particular location $\x$ can be chosen
as $I(\x) \propto f(\x)$. However, the attractiveness $\b$ is relative, it should be
seen in the eyes of the beholder or judged by the other fireflies. Thus, it will
vary with the distance $r_{ij}$ between firefly $i$ and firefly $j$.

In addition, light intensity decreases with the distance from its source, and light is also
absorbed in the media, so we should allow the attractiveness to vary with the degree of absorption. In the simplest form, the light intensity $I(r)$ varies according to the inverse
square law
\be I(r)=\frac{I_s}{r^2}, \ee where $I_s$ is the intensity at the source.
In a medium with a fixed light absorption coefficient $\gamma$,
the light intensity $I$ varies with the distance $r$. That is
\be I=I_0 e^{-\gamma r}, \ee
where $I_0$ is the original light intensity. In order to avoid the singularity
at $r=0$ in the expression $I_s/r^2$, the combined effect of both the inverse
square law and absorption can be approximated as the following Gaussian form
$I(r)=I_0 e^{-\gamma r^2}.$
As a firefly's attractiveness is proportional to
the light intensity seen by adjacent fireflies, we can now define
the attractiveness $\b$ of a firefly by
\be \b = \b_0 e^{-\gamma r^2}, \label{att-equ-100} \ee
where $\b_0$ is the attractiveness at $r=0$.
In fact, equation (\ref{att-equ-100}) defines a characteristic distance
$\Gamma=1/\sqrt{\gamma}$ over which the attractiveness changes significantly
from $\b_0$ to $\b_0 e^{-1}$.

The distance between any two fireflies $i$ and $j$ at $\x_i$ and $\x_j$, respectively, is
the Cartesian distance
\be r_{ij}=||\x_i-\x_j|| =\sqrt{\sum_{k=1}^d (x_{i,k} - x_{j,k})^2}, \ee
where $x_{i,k}$ is the $k$th component of the spatial coordinate $\x_i$ of $i$th
firefly. It is worth pointing out that the distance $r$ defined above is {\it not} limited
to the Euclidean distance. We can define other distance $r$ in the $n$-dimensional
hyperspace, depending on the type of problem of our interest. For example,
for job scheduling problems, $r$ can be defined as the time lag or time interval.
For complicated networks such as the Internet and social networks, the distance
$r$ can be defined as the combination of the degree of local clustering and
the average proximity of vertices. In fact, any measure that can effectively characterize
the quantities of interest in the optimization problem can be used as the `distance' $r$.

The movement of a firefly $i$ is attracted to another more attractive (brighter)
firefly $j$ is determined by
\be \x_i =\x_i + \b_0 e^{-\gamma r^2_{ij}} (\x_j-\x_i) + \a \; \ff{\e}_i,
\label{FA-equ-50}  \ee
where the second term is due to the attraction. The third term
is randomization with $\a$ being the randomization parameter, and
$\ff{\e}_i$ is a vector of random numbers drawn from a Gaussian distribution
or uniform distribution.  For most implementations,
we can take $\b_0=1$ and $\a=O(1)$.
It is worth pointing out that (\ref{FA-equ-50}) is a random walk biased
towards the brighter fireflies. If $\b_0=0$, it becomes a simple random walk.
Furthermore, the randomization term can easily be
extended to other distributions such as L\'evy flights\kkk{Gutow,Pav}.

The parameter $\gamma$ now characterizes the variation of the attractiveness,
and its value is crucially important in determining the speed of the convergence
and how the FA algorithm behaves. In theory, $\gamma \in [0,\infty)$, but in
practice, $\gamma=O(1)$ is determined by the characteristic length $\Gamma$ of the
system to be optimized. Thus, for most applications, it typically varies
from $10^{-5}$ to $10^{5}$.

\subsection{Asymptotic Limits}

The typical scale $\Gamma$ should be associated with the scale concerned in our
optimization problem. If $\Gamma$ is the typical scale for a given
optimization problem, for a very large number of fireflies $n \gg k$ where
$k$ is the number of local optima, then the initial locations of these $n$
fireflies should distribute relatively uniformly over
the entire search space. As the iterations proceed, the fireflies would converge
into all the local optima (including the global ones). By comparing the best
solutions among all these optima, the global optima can easily be achieved.

Two important limiting or asymptotic cases when $\gamma \rightarrow 0$ and $\gamma \rightarrow \infty$.
One limiting case is $\gamma \rightarrow 0$, the attractiveness is constant $\b=\b_0$ and $\Gamma \rightarrow \infty$,
this is equivalent to saying that the light intensity does not decrease in an idealized sky.
Thus, a flashing firefly can be seen anywhere in the domain. Thus, a single (usually global) optima can easily be reached. If we replace $\x_j$ by the current global best $\ff{g}_*$,
then the Firefly Algorithm becomes the special case of accelerated particle swarm optimization (PSO) discussed earlier.
Subsequently, the efficiency of this special case
is the same as that of PSO.

Another limiting case is $\gamma \rightarrow \infty$, leading to
$\Gamma \rightarrow 0$ and  $\b(r) \rightarrow \delta(r)$ which is the Dirac delta function,
which means that the attractiveness is almost zero in the sight of other fireflies. This is equivalent
to the case where the fireflies roam in a very thick foggy region randomly. No other fireflies can be
seen, and each firefly roams in a completely random way. Therefore, this corresponds
to the completely random search method.

In general, firefly algorithm usually works between these two extremes, and it is thus possible
to adjust the parameter $\gamma$ and $\alpha$ so that it can outperform both the random search
and PSO. In fact, FA can find the global optima as well as the local optima
simultaneously and effectively\kkk{Yang3}. A further advantage of FA is that different fireflies will work almost independently, it is
thus particular suitable for parallel implementation. It is even better than genetic algorithms
and PSO because fireflies aggregate more closely around each optimum.
It can be expected that  the interactions between different subregions are
minimal in parallel implementation.

\subsection{Convergence and Chaos}

We can carry out the convergence analysis for the firefly algorithm in a
framework similar to Clerc and Kennedy's dynamical analysis. For simplicity,
we start from the equation for firefly motion without the randomness term
\be \x_i^{t+1}=\x_i^t + \b_0 e^{-\gamma r_{ij}^2} (\x_j^t - \x_i^t). \ee
If we focus on a single agent, we can replace $\x_j^t$ by the global best
$g$ found so far, and we have
\be \x_i^{t+1} =\x_i^t + \b_0 e^{-\gamma r_{i}^2 } (g-\x_i^t), \ee
where the distance $r_i$ can be given by the $\ell_2$-norm
\be r_i^2=||g-\x_i^t||_2^2, \ee

\begin{figure}
\centerline{\includegraphics[height=2in,width=2.5in]{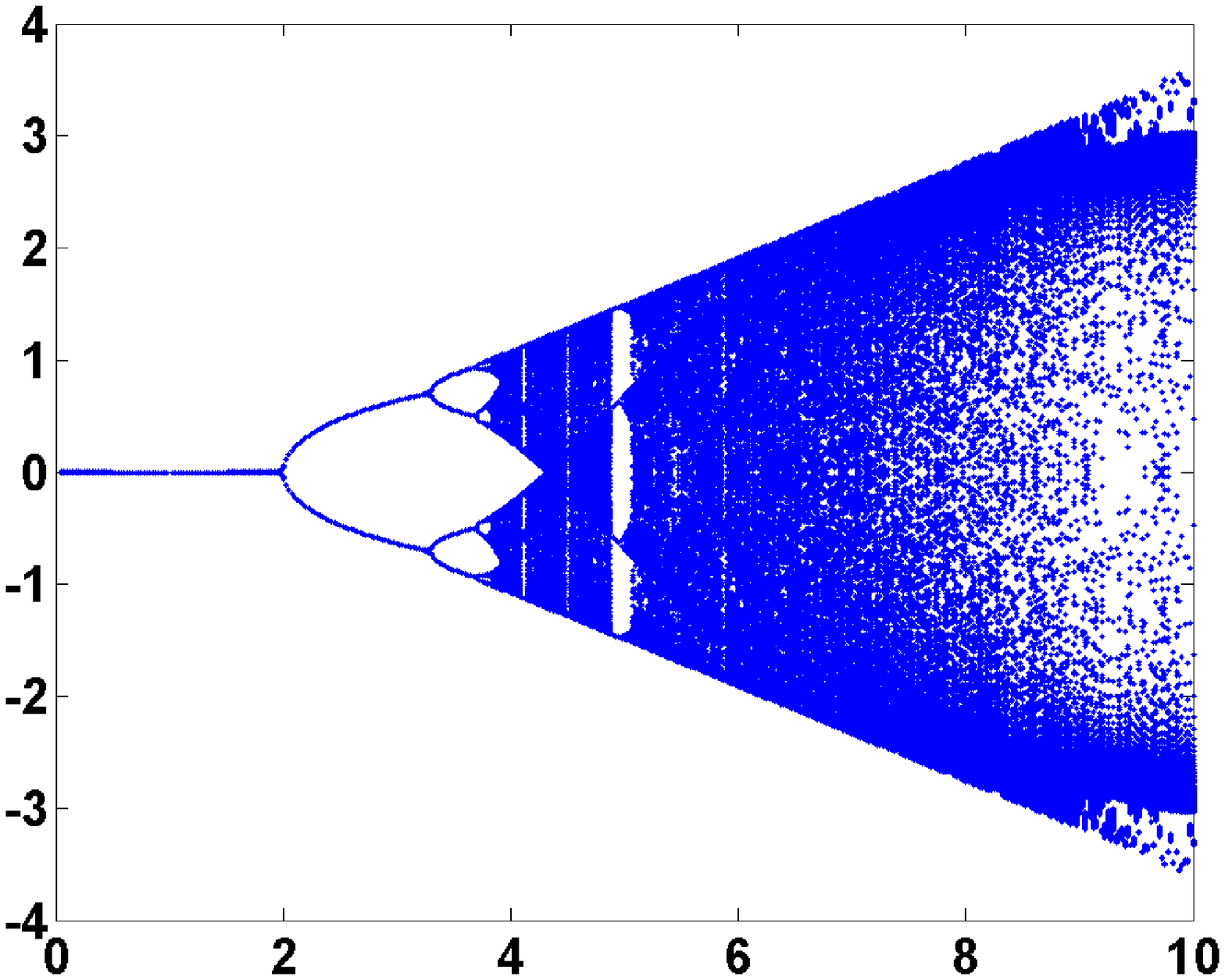}
\includegraphics[height=2in,width=3in]{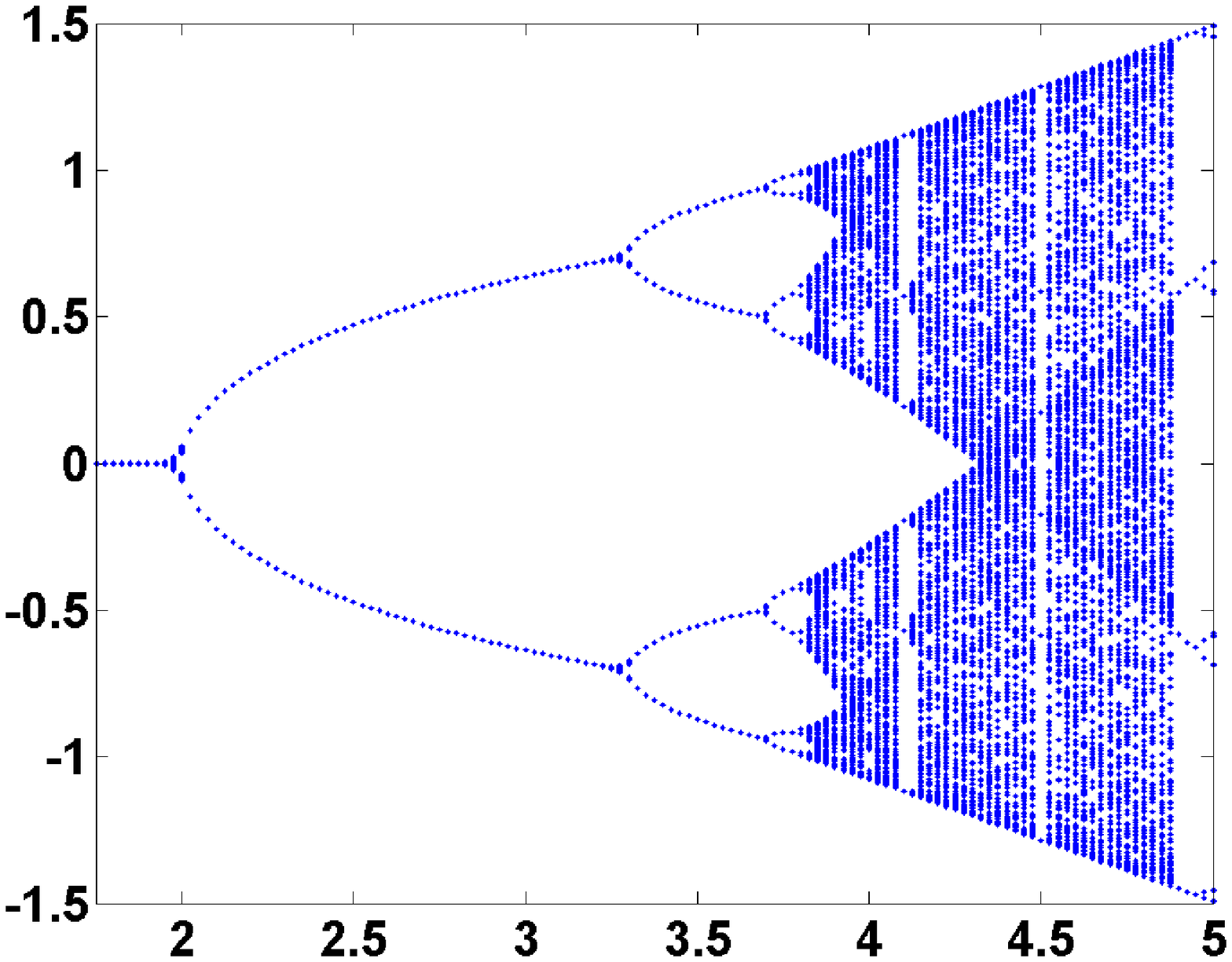}}
\caption{The chaotic map of the iteration formula (\ref{fa-dynamics-150}) in the firefly algorithm
and the transition between from periodic/multiple  states to chaos. \label{fap-fig-100} }
\end{figure}

In an even simpler 1-D case, we can set $y_{t}=g-\x_i^t$, and we have
\be y_{t+1} =y_t - \b_0 e^{-\gamma y_t^2} y_t. \label{fa-dynamics-25} \ee
Using the approximation of $0 < \gamma y_t^2 \ll 1$, we have
\be y_{t+1} = y_t -\b_0 (1-\gamma y_t^2) y_t = y_t [1-\b_0 (1-\gamma y_t^2)]. \label{fa-dynamics-50} \ee
We can see that $\gamma$ is a scaling parameter which only affects the scales/size
of the firefly movement. In fact, we can let
$u_t=\sqrt{\gamma} y_t$ and rewrite the above the equation as
\be u_{t+1} =u_t [1-\b_0 (1-u_t^2)], \label{fa-dynamics-100} \ee
or directly from (\ref{fa-dynamics-25})
\be u_{t+1}=u_t [1-\b_0 e^{-u_t^2}]. \label{fa-dynamics-150} \ee
These equations can be analyzed easily using the same methodology for studying
the well-known logistic map
\be u_t=\lam u_t (1-u_t). \label{fa-chaos-55} \ee
The chaotic map is shown
in Fig. \ref{fap-fig-100}, and the focus on the transition from periodic multiple
states to chaotic behaviour is shown in the same figure.

As we can see from Fig. \ref{fap-fig-100} that convergence can be achieved
for $\b_0 <2$. There is a transition from periodic to chaos at $\b_0 \approx 4$.
This may be surprising, as the aim of designing a metaheuristic algorithm is
to try to find the optimal solution efficiently and accurately. However,
chaotic behaviour is not necessarily a nuisance, in fact, we can use it to
the advantage of the firefly algorithm. Simple chaotic characteristics from  (\ref{fa-chaos-55})
can often be used as an efficient mixing technique for generating diverse solutions.
Statistically, the logistic mapping (\ref{fa-chaos-55}) for $\lam=4$ for the initial
states in (0,1) corresponds a beta distribution
\be B(u,p,q)=\frac{\Gamma(p+q)}{\Gamma(p) \Gamma(q)} u^{p-1} (1-u)^{q-1}, \ee
when $p=q=1/2$. Here $\Gamma(z)$ is the Gamma function
\be \Gamma(z) =\int_0^{\infty} t^{z-1} e^{-t} dt.  \ee
In the case when  $z=n$ is an integer, we have $\Gamma(n)=(n-1)!$. In addition,
$\Gamma(1/2)=\sqrt{\pi}$. From the algorithm implementation point of view,
we can use higher attractiveness $\b_0$ during the early stage
of iterations so that the fireflies can explore, even chaotically,
the search space more effectively. As the search continues and convergence
approaches, we can reduce the attractiveness $\b_0$ gradually, which may
increase the overall efficiency of the algorithm. Obviously, more
studies are highly needed to confirm this.

\section{Essential Components of Metaheuristics}

The efficiency of metaheuristic algorithms can be attributed to the fact
that they imitate the best features in nature, especially the selection of
the fittest in biological systems which have evolved
by natural selection over millions of years.

Metaheuristics can be considered as an efficient way to produce acceptable
solutions by trial and error to a complex problem in
a reasonably practical time. The complexity of the problem of interest makes it impossible
to search every possible solution or combination, the aim is
to find good feasible solution in an acceptable timescale. There is
no guarantee that the best solutions can be found, and we even do not
know whether an algorithm will work and why if it does work. The idea is
to have an efficient but practical algorithm that will work most
the time and is able to produce good quality solutions. Among the
found quality solutions, it is expected some of them are nearly optimal,
though there is no guarantee for such optimality.

The main components of any metaheuristic algorithms are: intensification
and diversification, or exploitation and exploration\kkk{Blum,Yang2}.
Diversification means to generate diverse solutions so as to explore the search space
on the global scale, while intensification means to focus on the
search in a local region by exploiting the information that a current
good solution is found in this region. This is in combination with
with the selection of the best solutions.

The selection of the best ensures
that the solutions will converge to the optimality.
On the other hand, the diversification
via randomization avoids the solutions being trapped at local optima, while increases the diversity
of the solutions. The good combination of these two major components will usually ensure
that the global optimality is achievable.

The fine balance between these two components is very
important to the overall efficiency and performance of an algorithm.
Too little exploration and too much exploitation could cause the
system to be trapped in local optima, which makes it very difficult
or even impossible to find the global optimum. On the other hand,
if too much exploration but too little exploitation, it may be difficult
for the system to converge and thus slows down the overall search performance.
The proper balance itself is an optimization problem,
and one of the main tasks of designing new algorithms is
to find a certain balance concerning this optimality and/or tradeoff.
Furthermore, just exploitation and exploration are not enough.
During the search, we have to use a proper mechanism or
criterion to select the best solutions. The  most common criterion
is to use the survival of the fittest, that is to keep updating the
the current best found so far. In addition, certain elitism is often
used, and this is to ensure the best or fittest solutions are not
lost, and should be passed onto the next generations.

\section{Randomization Techniques}

As discussed earlier, an important component in swarm intelligence
and modern metaheuristics is randomization, which enables an algorithm
to have the ability to jump out of any local optimum so as to search globally.
Randomization can also be used for local search around the current best
if steps are limited to a local region. Fine-tuning the randomness
and balance of local search and global search is crucially important
in controlling the performance of any metaheuristic algorithm.

Randomization techniques can be a very simple method using uniform distributions,
or more complex methods as those used in Monte Carlo simulations\kkk{Sobol}. They can also
be more elaborate, from Brownian random walks to L\'evy flights.

\subsection{Random Walks}

A random walk is a random process which consists of taking a series of consecutive
random steps. Mathematically speaking, let $S_N$ denotes the sum of each
consecutive random step $X_i$, then $S_N$ forms a random walk
\be S_N=\sum_{i=1}^N X_i = X_1 + ... + X_N, \ee
where $X_i$ is a random step drawn from a random distribution. This relationship can also be
written as a recursive formula
\be S_N=\sum_{i=1}^{N-1} + X_N = S_{N-1} + X_N, \label{Walk-markov-50}  \ee
which means the next state $S_N$ will only depend the current existing state $S_{N-1}$
and the motion or transition $X_N$ from the existing state to the next state.
This is typically the main property of a Markov chain to be introduced later.
Here the step size or length in a random walk can be fixed or varying.
Random walks have  many applications in physics, economics, statistics,
computer sciences, environmental science and engineering.

In theory, as the number of steps $N$ increases, the central limit theorem implies that
the random walk (\ref{Walk-markov-50}) should approaches a Gaussian distribution.
In addition, there is no reason why each step length
should be fixed. In fact, the step size can also vary according to a known distribution.
If the step length obeys the Gaussian distribution, the random walk becomes the
Brownian motion\kkk{Gutow,Yang2}.

As the mean of particle locations
is obviously zero, their variance will increase linearly with $t$.
In general, in the $d$-dimensional space, the variance of Brownian random walks can be written as
\be \sigma^2(t) = |v_0|^2 t^2 + (2 d D) t, \ee
where $v_0$ is the drift velocity of the system. Here $D=s^2/(2 \tau)$ is the effective diffusion coefficient
which is related to the step length $s$ over a short time interval $\tau$ during each jump.

Therefore, the Brownian motion $B(t)$ essentially obeys a Gaussian distribution with zero mean and time-dependent variance.
That is, $ B(t) \sim N(0, \sigma^2(t))$
where $\sim$ means the random variance obeys the distribution on the right-hand side; that is, samples should be drawn
from the distribution. A diffusion process can be viewed as a series of Brownian motion, and the motion
obeys the Gaussian distribution. For this reason, standard diffusion is often referred to as the Gaussian diffusion.
If the motion at each step is not Gaussian, then the diffusion is called non-Gaussian diffusion.
If the step length obeys other distribution, we have to deal with more generalized random walk. A very special case
is when the step length obeys the L\'evy distribution, such a random walk is called  L\'evy flight or L\'evy walk.

\subsection{L\'evy Distribution and L\'evy Flights}

In nature, animals search for food in a random or quasi-random manner.
In general, the foraging path of an animal is effectively a random walk because the next move is
based on the current location/state and the transition probability to
the next location. Which direction it chooses depends implicitly on a
probability which can be modelled mathematically. For example, various studies
have shown that the flight behaviour of many animals and insects has demonstrated
the typical characteristics of L\'evy flights\kkk{Pav,Romos,Reynolds,Reynolds2}.
For example, a recent study by Reynolds and Frye shows that
fruit flies or {\it Drosophila melanogaster}, explore their landscape using a series
of straight flight paths punctuated by a sudden $90^{o}$ turn, leading to
a L\'evy-flight-style intermittent scale-free search pattern\kkk{Reynolds}.
Subsequently, such behaviour has been applied to
optimization and optimal search, and preliminary results show its
promising capability \kkk{Viswan,Pav}.

Broadly speaking, L\'evy flights are a random walk whose step length is drawn from the L\'evy distribution,
often in terms of a simple power-law formula $L(s) \sim |s|^{-1-\b} $ where $0<\b \le 2$ is an index.
Mathematically speaking, a simple version of L\'evy distribution can be defined as
\be L(s, \gamma, \mu)=\brak{\sqrt{\frac{\gamma}{2 \pi}} \exp[-\frac{\gamma}{2 (s -\mu)}] \frac{1}{(s-\mu)^{3/2}},
& 0< \mu <  s < \infty \\ \\ 0 & \textrm{otherwise},} \ee
where $\mu>0$ is a minimum step and $\gamma$ is a scale parameter. Clearly, as $s \rightarrow \infty$,
we have \be L(s, \gamma, \mu) \= \sqrt{\frac{\gamma}{2 \pi}} \frac{1}{s^{3/2}}. \ee
This is a special case of the generalized L\'evy distribution\kkk{Gutow,Nol}.

In general, L\'evy distribution should be defined in terms of Fourier transform
\be F(k)=\exp[-\a |k|^{\b}], \quad 0 < \b \le 2, \ee
where $\a$ is a scale parameter. The inverse of this integral is not easy, as it does not
have analytical form, except for a few special cases.

For the case of $\b=2$, we have \be F(k)=\exp[-\a k^2], \ee
whose inverse Fourier transform corresponds to a Gaussian distribution.
Another special case is $\b=1$, and we have
\be F(k)=\exp[-\a |k|], \ee
which corresponds to a Cauchy distribution
\be p(x,\gamma,\mu)=\frac{1}{\pi} \frac{\gamma}{\gamma^2+(x-\mu)^2}, \ee
where $\mu$ is the location parameter, while $\gamma$ controls the
scale of this distribution.

For the general case, the inverse integral
\be L(s) =\frac{1}{\pi} \int_0^{\infty} \cos (k s) \exp[-\a |k|^{\b}] dk, \ee
can be estimated only when $s$ is large. We have
\be L(s) \rightarrow \frac{\a \; \b \; \Gamma(\b) \sin (\pi \b/2)}{ \pi |s|^{1+\b}}, \quad s \rightarrow \infty. \ee

L\'evy flights are more efficient than Brownian random walks
in exploring unknown, large-scale search space.
There are many reasons to explain this efficiency, and one of them
is due to the fact that the variance of L\'evy flights
\be \sigma^2(t) \sim t^{3-\b},  \quad 1 \le \b \le 2, \ee
increases much faster than the linear relationship (i.e., $\sigma^2(t) \sim t$) of Brownian random walks.

From the implementation point of view, the generation of random numbers with L\'evy flights
consists of two steps: the choice of a random direction and the generation of
steps which obey the chosen  L\'evy distribution. The generation of a direction should be
drawn from a uniform distribution, while the generation  of steps is quite tricky.
There are a few ways of achieving this, but one of the most efficient and yet straightforward
ways is to use the so-called Mantegna algorithm for a symmetric L\'evy stable distribution\kkk{Mant,Nol}.
Here `symmetric' means that the steps can be positive and negative.

A random variable $U$ and its probability distribution can be called stable if
a linear combination of its two identical copies (or $U_1$ and $U_2$) obeys the
same distribution. That is, $a U_1 + b U_2$ has the same distribution as
$c U + d$ where $a,b>0$ and $c,d \in \Re$. If $d=0$, it is called strictly stable.
Gaussian, Cauchy and L\'evy distributions are all stable distributions.

\begin{figure}
\centerline{\includegraphics[height=1.5in,width=4in]{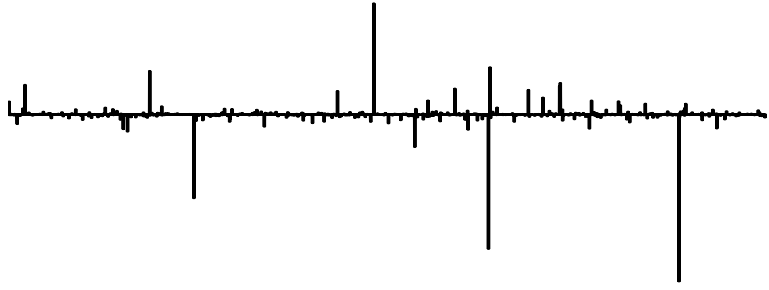}}
\caption{Steps drawn from a L\'evy distribution ($\b=1$) using
Mantegna's algorithm.  \label{fap-fig-300} }
\end{figure}

In Mantegna's algorithm\kkk{Mant}, the step length $s$ can be calculated by
\be s=\frac{u}{|v|^{1/\b}}, \ee
where $u$ and $v$ are drawn from normal distributions. That is
\be u \sim N(0, \sigma_u^2), \quad v \sim N(0, \sigma_v^2), \ee
where
\be \sigma_u=\Big\{\frac{\Gamma(1+\b) \sin (\pi \b/2)}{\Gamma[(1+\b)/2] \; \b \; 2^{(\b-1)/2} }\Big\}^{1/\b},
\quad \sigma_v=1. \ee
This distribution (for $s$) obeys the expected L\'evy distribution for $|s| \ge |s_0|$ where
$s_0$ is the smallest step. In principle, $|s_0| \gg 0$, but in reality $s_0$ can be
taken as a sensible value such as $s_0=0.1$ to $1$.

\begin{figure}
\centerline{\includegraphics[height=2in,width=4in]{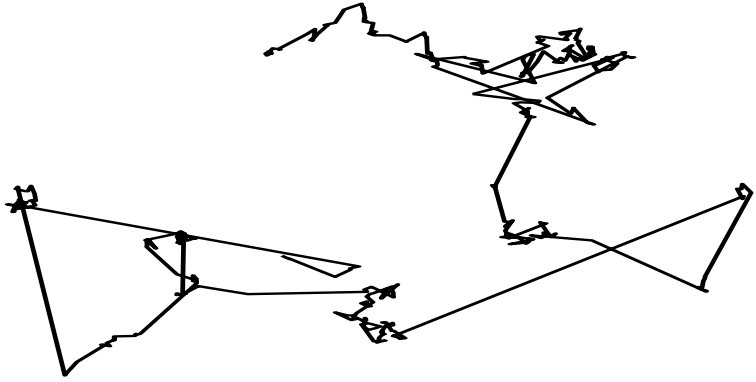}}
\caption{Schematic representation  of 2D L\'evy flights.  \label{fap-fig-400} }
\end{figure}

Fig. \ref{fap-fig-300} shows the 250 steps drawn from a L\'evy distribution
with $\b=1$, while Fig. \ref{fap-fig-400} shows the path of L\'evy flights of 250 steps starting from $(0,0)$.
It is worth pointing out that a power-law distribution is often linked to
some scale-free characteristics, and L\'evy flights can thus show
self-similarity and fractal behavior in the flight patterns.

Studies show that L\'evy flights can maximize the efficiency of resource searches in uncertain environments.
In fact,  L\'evy flights have been observed among
foraging patterns of albatrosses and fruit flies, and spider monkeys\kkk{Pav,Romos,Reynolds,Reynolds2,Viswan}.  In addition,
L\'evy flights have many applications. Many physical phenomena such as
the diffusion of fluorenscent molecules, cooling behavior and noise could show L\'evy-flight
characteristics under the right conditions\kkk{Reynolds2}.

\subsection{Step Size in Random Walks \label{RW-step-555}}

As random walks are widely used for randomization and local search,
a proper step size is very important. In the generic equation
\be \x^{t+1}=\x^t + s \; \ff{\e}_t, \ee
$\ff{\e}_t$ is drawn from a standard normal distribution with
zero mean and unity standard deviation. Here the step size $s$
determines how far a random walker (e.g., an agent or particle
in metaheuristics) can go for a fixed number of iterations.

If $s$ is too large, then the new solution $\x^{t+1}$ generated will be
too far away from the old solution (or more often the current best).
Then, such a move is unlikely to be accepted. If $s$ is too small,
the change is too small to be significant, and consequently such
search is not efficient. So a proper step size is important to maintain
the search as efficient as possible.

From the theory of simple isotropic random walks, we know that
the average distance $r$ traveled in the $d$-dimension space is
\be r^2=2 d D t, \ee
where $D=s^2/2\tau$ is the effective diffusion coefficient.
Here $s$ is the step size or distance traveled at each jump,
and $\tau$ is the time taken for each jump. The above equation
implies that
\be s^2=\frac{\tau \; r^2}{t \; d}. \ee
For a typical length scale $L$ of a dimension of interest, the local search
is typically limited in a region of $L/10$. That is, $r=L/10$.
As the iterations are discrete, we can take $\tau=1$.
Typically in metaheuristics, we can expect that the number
of generations is usually $t=100$ to $1000$, which means that
\be s \=\frac{r}{\sqrt{t d}}=\frac{L/10}{\sqrt{t \; d}}. \label{stepsize-equ-555} \ee
For $d=1$ and $t=100$, we have $s=0.01L$, while $s=0.001L$ for $d=10$ and $t=1000$.
As step sizes could differ from variable to variable, a step size
ratio $s/L$ is more generic. Therefore, we can use $s/L=0.001$ to $0.01$ for most problems.

\subsection{Swarm Intelligence and Markov Chains}

From the discussion and analysis of metaheuristic algorithms, we know that
there is no mathematical framework in general to provide insights into the
working mechanisms, the stability and convergence of a give algorithm.
Despite the increasing popularity of metaheuristic, mathematical analysis
remains fragmental, and many open problems need urgent attention.

From the statistical point of view, most metaheuristic algorithms can be viewed
in the framework of Markov chains\kkk{Fishman,Sobol}. For example,
simulated annealing\kkk{Kirk} is a Markov
chain, as the next state or new solution
in SA only depends on the current state/solution and the transition probability.
For a given Markov chain with certain ergodicity, a stability probability distribution
and convergence can be achieved. In fact, it has been proved that SA
will always converge if the system is cooled down slowly enough and the
annealling process runs long enough.

Now if look at the PSO closely using the framework of Markov chain Monte Carlo\kkk{Gamer,Geyer,Ghate}, each particle in PSO essentially forms a Markov chain, though this Markov chain is biased towards to the current best, as the transition probability often leads to the
acceptance of the move towards the current global best. In addition, the multiple
Markov chains are interacting in terms of partly deterministic attraction movement.
Therefore, the mathematical analysis concerning of the rate of convergence of
PSO is very difficult, if not impossible. Similarly, other swarm-based metaheuristics
such as firefly algorithms and cuckoo search can also be viewed as interacting
Markov chains. There is no doubt that
any theoretical advance in understanding multiple interacting Markov chains
will gain tremendous insight  in understanding how the swarm intelligence
behaves and may consequently lead to the design of better or new metaheuristics.

\section{Discussions}

Swarm-based metaheuristics are becoming increasingly popular, though theoretical
framework and analysis still lack behind. Many important problems still remain
unsolved. Firstly, there is no general convergence analysis for all metaheuristics,
apart from some fragmental and yet important work concerning PSO and SA.
Secondly, the choice of algorithm-dependent parameters such as population size
and probability is currently based on experience and parametric studies.
Thirdly, stopping criteria are also a bit arbitrary, though fixed accuracy
or maximum number of iterations are widely used.

In addition, there is no agreed measure to compare the performance
of any two algorithms. People either compare the number of function
evaluations for a given objective value or compare the objective values
for a given number of functional evaluations. These results are often
affected by algorithm-dependent parameters or way of actual
implementation. As randomness is an integrated part of metaheuristics,
statistical measures should be used. A general, statistical framework
for performance comparison and algorithm evaluations is highly needed.

Furthermore, no-free-lunch theorems\kkk{Wolpert} may have some important
implications in metaheuristics, however, the basic assumptions of these
theorems may be questionable. Besides, these theorems are valid for
single objective optimization, and they remain open questions for multiobjective.

Finally, the current understanding of swarm intelligence is still very limited.
How multiagents interact with simple rules can lead to complex behaviour,
self-organization and eventually intelligence still remains
a major mystery. Any insight with theoretical basis may have important
implications to many areas, not only metaheuristics and optimization,
but also neural networks and machine learning. Whatever the future developments
will be, there is no doubt that swarm-based algorithms and swarm intelligence
will play an ever more important role in science and engineering.

\end{document}